\documentclass[review]{elsarticle}
\usepackage{amsthm}
\usepackage{amsmath}
\usepackage{color}
\usepackage{lineno,hyperref}
\newtheorem{proposition}{Proposition}[section]
\newtheorem{lemma}[proposition]{Lemma}

\newtheorem{definition}[proposition]{Definition}

\newtheorem{theorem}[proposition]{Theorem}

\newtheorem{remark}[proposition]{Remark}
\newcommand{\R}{\mathbf{R}}

\newcommand{\ra}{{\rightarrow}}

\modulolinenumbers[5]

\journal{Journal of \LaTeX\ Templates}

\begin{document}

\begin{frontmatter}

\title{The existence of the solution of the wave equation on graphs
}


\author[mymainaddress]{Yong Lin}
\ead{yonglin@tsinghua.edu.cn}

\author[mysecondaryaddress]{Yuanyuan Xie\corref{mycorrespondingauthor}}
\cortext[mycorrespondingauthor]{Corresponding author}
\ead{yyxiemath@163.com}

\address[mymainaddress]{Yau Mathematical Sciences Center, Tsinghua University, Beijing, 100084, China.}
\address[mysecondaryaddress]{School of Mathematics, Renmin University of China, Beijing, 100872, China.}

\begin{abstract}

Let $G=(V, E)$ be a finite weighted graph, and $\Omega\subseteq V$ be a domain such that $\Omega^\circ\neq\emptyset$. In this paper, we study the following initial boundary problem for the non-homogenous wave equation
\begin{equation*}
\left\{
\begin{aligned}
&\partial_t^2 u(t,x)-\Delta_\Omega u(t,x)=f(t,x),\qquad&&(t,x)\in[0,\infty)\times \Omega^\circ,\\
&u(0,x)=g(x),\qquad&& x\in\Omega^\circ,\\
&\partial_tu(0,x)=h(x),\qquad&& x\in\Omega^\circ,\\
&u(t,x)=0,\qquad&&(t,x)\in[0,\infty)\times\partial \Omega,
\end{aligned}
\right.
\end{equation*}
where $\Delta_\Omega$ denotes the Dirichlet Laplacian on $\Omega^\circ$. Using Rothe's method, we prove that the above wave equation has a unique solution.

\end{abstract}


\end{frontmatter}

\linenumbers

\section{Introduction}

\setcounter{equation}{0}

The wave equation on graphs was first studied by Courant et al.\cite{Courant-Friedrichs-Lewy_1928}.
On metric graphs, Friedman and Tillich \cite{Friedman-Tillich_2004} developed a wave equation that is based on the edge-based Laplacian, and also gave some applications of wave equation; Schrader \cite{Schrader_2009} studied the solution of the wave equation, and established the finite propagation speed.

In recent years, the study of equations on graphs has attracted attention of researchers in various fields.
Grigoryan et al. \cite{Grigoryan-Lin-Yang_2016_1, Grigoryan-Lin-Yang_2017}, using the mountain pass theorem due to Ambrosetti-Rabinowitz, studied the existence of solutions to Yamabe type equation and some nonlinear equations on graphs, respectively. They \cite{Grigoryan-Lin-Yang_2016_2} also considered the Kazdan-Warner equation on graph. The proof uses the calculus of variations and a method of upper and lower solutions. The Kazdan-Warner equation has also been studied by Keller and Schwarz on canonically compactifiable graphs \cite{Keller-Schwarz_2018}. In \cite{Lin-Wu_2017}, Lin and Wu proved the existence and nonexistence of global solutions of the Cauchy Problem for $\partial_tu=\Delta u+u^{1+\alpha}$ with $\alpha>0$ on a finite or locally finite graph.

The purpose of this paper is to study the wave equation on graphs, which is defined by the Dirichlet Laplacian
(see \eqref{e:Delta_Omega_2}). In contrast to \cite{Friedman-Tillich_2004}, our Laplacian is defined on vertices.
It is well that variational method is commonly used to study the solution of elliptic boundary value problems, and based on the minimalization of the corresponding functional. However, this method is not used in parabolic and hyperbolic equations.
In 1930, Rothe's method was originally introduced by Rothe for the study of parabolic equations \cite{Rothe_1930}. After then, this method was used by many authors for parabolic and hyperbolic problems, see for example Rektorys \cite{Rektorys_1971}, Kacur \cite{Kacur_1984} and so on. In this paper, using Rothe's method, we prove that the non-homogeneous wave equation has a unique solution on finite weighted graphs.

Let $G=(V, E)$ be a finite weighted graph, where $V$ and $E$ denote the vertex set and the edge set of $G$, respectively. Given a non-empty domain (or connected set) $\Omega\subseteq V$, the boundary of $\Omega$ is defined by
\begin{eqnarray}\label{e:boundary_Omega}
\partial \Omega:=\{x\in\Omega: \mbox{ there exists } y\in \Omega^c \mbox{ such that } y\sim x \},
\end{eqnarray}
and the interior of $\Omega$ is defined by $\Omega^\circ:=\Omega\backslash\partial\Omega$.
Assume that $\Omega^\circ\neq\emptyset$.
We consider the following initial boundary value problem
\begin{equation}\label{e:wave_eq_on_graph}
\left\{
\begin{aligned}
&\partial_t^2u(t,x)-\Delta_\Omega u(t,x)=f(t,x),&\quad(t,x)\in[0,\infty)\times \Omega^\circ,\\
&u(0,x)=g(x),&\quad x\in\Omega^\circ,\\
&\partial_tu(0,x)=h(x),&\quad x\in\Omega^\circ,\\
&u(t,x)=0,&\quad(t,x)\in[0,\infty)\times\partial \Omega,
\end{aligned}
\right.
\end{equation}
where $\Delta_\Omega$ is the Dirichlet Laplacian on $\Omega^\circ$, $u:[0,\infty)\times\Omega^\circ\to \R$ is unknown,
$f:[0,\infty)\times\Omega^\circ\to\R$ is continuous with respect to $t$, and $g, h: \Omega^\circ\to \R$.

\begin{definition}
A function $u:[0,\infty)\times\Omega\to\R$ is said to be a {\em solution} of \eqref{e:wave_eq_on_graph} on $[0,\infty)\times\Omega$ if $u$ is twice continuously differentiable with respect to $t$, and \eqref{e:wave_eq_on_graph} is satisfied.
\end{definition}

In the following, we briefly introduce Rothe's method. Details appears in Section \ref{S:proof_of_main_results}.
For any $T>0$, we divide the time interval $[0,T]$ into $n$ subintervals $[t_{i-1}, t_i]$ that have the same length $\ell=T/n$, where $t_i=i\ell$ for $i=0,\ldots, n$.
Let $u_{n,i}(x)$ be the approximation of a solution $u(t,x)$ of \eqref{e:wave_eq_on_graph} at $t=t_i$. Then we can replace $\partial_t u(t_i,x)$ and $\partial_t^2 u(t_i,x)$ by the difference quotients
$$
\delta u_{n,i}=(u_{n,i-1}-u_{n,i})/{\ell}\quad\mbox{and}
\quad\delta^2 u_{n,i}=(\delta u_{n,i}-\delta u_{n,i-1})/{\ell},$$
respectively. Using the functions $\{u_{n,i}\}_{i=1}^n$, we construct the Rothe's functions $u^{(n)}(t,x)$ (see \eqref{e:defi_u(n)}) and some auxiliary functions. Then we show that $u^{(n)}(t,x)$ converges to a solution $u(t,x)$ of \eqref{e:wave_eq_on_graph}.

Our main results are as follows:

\begin{theorem}\label{T:sol_of_wave_equ}
Let $G=(V, E)$ be a finited weighted graph, $\Omega\subseteq V$ be a domain such that $\Omega^\circ\neq \emptyset$, and $\Delta_\Omega$ be the Dirichlet Laplacian on $\Omega^\circ$. Assume that there exist some positive constants $\alpha$ and $c:=c(\Omega^\circ)$ such that
\begin{equation}\label{eq:condition_on_f}
\|f(t,\cdot)-f(s,\cdot)\|_{L^2(\Omega^\circ)}\le c|t-s|^\alpha\qquad \mbox{ for any }t,s\in[0,\infty).
\end{equation}
Then \eqref{e:wave_eq_on_graph} has a unique solution.
\end{theorem}

\begin{remark}
\begin{enumerate}
\item[(1)] The assumption \eqref{eq:condition_on_f} is used to prove the existence of a solution of \eqref{e:wave_eq_on_graph}.
\item[(2)] It is easy to see that $f(t,x)=t\cdot\phi(x)$ and $f(t,x)=\phi(x)\cdot\sin t$ satisfy the hypotheses of
     Theorem \ref{T:sol_of_wave_equ}, where $\phi: V \rightarrow \R$.
\end{enumerate}
\end{remark}

\begin{theorem}\label{T:sol_of_nonlinear_wave_equ}
Let $G=(V, E)$ be a finite weighted graph, $\Omega\subseteq V$ be a domain such that $\Omega^\circ\neq \emptyset$, and $\Delta_\Omega$ be the Dirichlet Laplacian on $\Omega^\circ$. Also let $\{\varphi_k\}_{k=1}^N$ be an orthonormal basis of $W_0^{1,2}(\Omega)$ consisting of the eigenfunctions of $-\Delta_\Omega$ such that $-\Delta_\Omega\varphi_k=\lambda_k\varphi_k$ for $k=1,\ldots, N$, where $N=\#\Omega^\circ$. If
\begin{eqnarray}\label{eq:form_of_fgh}
\quad f(t,x)=\sum_{k=1}^N b_k(t)\varphi_k(x),\quad g(x)=\sum_{k=1}^Ng_k\varphi_k(x),\quad h(x)=\sum_{k=1}^Nh_k\varphi_k(x),
\end{eqnarray}
then the solution of \eqref{e:wave_eq_on_graph} is given by
\begin{eqnarray*}
\begin{aligned}
u(t,x)
=&\sum_{k=1}^N\frac{1}{\sqrt{\lambda_k}}\int_0^t\sin(\sqrt{\lambda_k}(t-s))b_k(s)\,ds\cdot \varphi_k(x)\\
&+\sum_{k=1}^Ng_k\cos(\sqrt{\lambda_k}t)\varphi_k(x)
+\sum_{k=1}^N\frac{1}{\sqrt{\lambda_k}}h_k\sin(\sqrt{\lambda_k}t)\varphi_k(x)
\end{aligned}
\end{eqnarray*}
\end{theorem}

The rest of the paper is organized as follows. In Section \ref{S:pre}, we introduce some definitions and
notations. Section \ref{S:proof_of_main_results} is devoted to the proof of Theorem \ref{T:sol_of_wave_equ}. In Section \ref{S:proof_of_main_results_nonlinear}, we give the proof of Theorem~\ref{T:sol_of_nonlinear_wave_equ}.

\section{Preliminaries}\label{S:pre}
\setcounter{equation}{0}

Let $G=(V, E)$ be a finite graph. We write $y\sim x$ if $xy\in E$. For any edge $xy\in E$, we assume that its weight $\omega_{xy}> 0$ and $\omega_{xy}=\omega_{yx}$. A pair $(V,\omega)$ is called a \textit{weighted graph}.
Furthermore, let $\mu: V\to \R^+$ be a positive finite measure. Define
\begin{eqnarray*}\label{e:D_mu}
D_\mu:=\max_{x\in V}\frac{m(x)}{\mu(x)},
\end{eqnarray*}
where $m(x):=\sum\limits_{y\sim x}\omega_{xy}$.
$G$ is called \textit{connected} if for any $x, y\in V$, there is a path connecting them.
In this paper, we consider finite weighted graphs.

Let $C(V)$ be the set of real functions on $V$. For any $u\in C(V)$, the \textit{$\mu$-Laplacian} $\Delta$ of $u$ is defined as follows:
\begin{equation*}\label{e:mu_Laplacian}
\Delta u(x)=\frac{1}{\mu(x)}\sum_{y\sim x}\omega_{xy}\big(u(y)-u(x)\big).
\end{equation*}
The associated \textit{gradient form} is defined by
\begin{equation*}\label{e:Gamma}
\Gamma(u,v)(x)=\frac{1}{2\mu(x)}\sum_{y\sim x}\omega_{xy}\big(u(y)-u(x)\big)\big(v(y)-v(x)\big).
\end{equation*}
Write $\Gamma(u)=\Gamma(u,u)$.
We denote the \textit{length of its gradient} by
\begin{equation*}\label{e:length_of_gradient}
|\nabla u|(x)=\sqrt{\Gamma(u)(x)}=\Big(\frac{1}{2\mu(x)}\sum_{y\sim x} \omega_{xy}\big(u(y)-u(x)\big)^2\Big)^{1/2}.
\end{equation*}

An integral of $u\in C(V)$ is defined by
\begin{equation*}\label{e:integration_of_function}
\int_V u\,d\mu=\sum_{x\in V}u(x)\mu(x).
\end{equation*}

Given a non-empty domain $\Omega\subseteq V$, let $\partial\Omega$ be defined as in \eqref{e:boundary_Omega}, and $\Omega^\circ=\Omega\backslash{\partial\Omega}$. For any $u\in C(\Omega^\circ)$, the \textit{Dirichlet Laplacian} $\Delta_\Omega$ on $\Omega^\circ$ is defined as follows: first we extend $u$ to the whole $V$ by setting $u\equiv 0$ outside $\Omega^\circ$ and then set
\begin{eqnarray*}\label{e:Delta_Omega_1}
\Delta_\Omega u=(\Delta u)|_{\Omega^\circ}.
\end{eqnarray*}
Then
\begin{eqnarray}\label{e:Delta_Omega_2}
\Delta_\Omega u(x)=\frac{1}{\mu(x)}\sum_{y\sim x}\omega_{xy}\big(u(y)-u(x)\big)\quad\mbox{ for any } x\in \Omega^\circ,
\end{eqnarray}
where $u(y)=0$ whenever $y\notin \Omega^\circ$. It is easy to see that $-\Delta_\Omega$ is a positive self-adjoint operator (see \cite{Grigoryan_2009, Weber_2012}).

\begin{lemma}\label{L:Green_formula}(Green's formula)\cite{Grigoryan_2009}
For any $u, v\in C(\Omega^\circ)$, we have
\begin{eqnarray*}
\int_{\Omega^\circ}\Delta_\Omega u\cdot v\,d\mu=-\int_\Omega\Gamma(u,v)\,d\mu.
\end{eqnarray*}
\end{lemma}

For given $\Omega\subseteq V$, let $L^2(\Omega)$ and the norm be
$$L^2(\Omega):=\{u\in C(V):\int_\Omega|u|^2\,d\mu<\infty\}\qquad\mbox{and}\qquad
\|u\|^2_{L^2(\Omega)}:=\int_\Omega|u|^2\,d\mu,$$
respectively. Also let $C_0(\Omega)$ be a set of all functions $u:\Omega\to \mathbf{R}$ with $u|_{\partial\Omega}=0$.

\section{Proof of Theorem~\ref{T:sol_of_wave_equ}}\label{S:proof_of_main_results}
\setcounter{equation}{0}

In this section, we apply Rothe's and energy methods to prove Theorem \ref{T:sol_of_wave_equ}, which states the existence and uniqueness of the solution for problem \eqref{e:wave_eq_on_graph}.

\subsection{Some priori estimates}\label{SS:estimates}

In this subsection, we give some lemmas that will be used to prove the existence of a solution.
For any $T>0$, choose a fixed integer $n>0$, we divide the interval $[0,T]$ into $n$ subintervals $[t_{i-1}, t_i]$ of same length $\ell$, where $t_0=0$, $t_i=i\ell$ for $i=1,\cdots, n$ and $n\ell=T$. Setting
$$
u_{n,0}(x):=
\left\{
\begin{aligned}
&g(x), \quad&&\mbox{ if }x\in\Omega^\circ,\\
&0,\qquad&&\mbox{ if }x\in\partial\Omega.
\end{aligned}
\right.
\quad
u_{n,-1}(x):=
\left\{
\begin{aligned}
&g(x)-\ell\cdot h(x),\quad&&\mbox{ if }x\in\Omega^\circ,\\
&0,\quad\qquad&&\mbox{ if }x\in\partial\Omega.
\end{aligned}
\right.
$$
and for $i=1,\ldots, n$,
$$
f_{n,i}(x):=
\left\{
\begin{aligned}
&f(t_i,x), \quad&&\mbox{ if }x\in\Omega^\circ,\\
&0,\quad&&\mbox{ if }x\in\partial\Omega.
\end{aligned}
\right. $$
Consider
\begin{eqnarray}\label{e:BVP_1}
(u-2u_{n,0}+u_{n,-1})/{\ell^2}-\Delta_\Omega u=f_{n,1}.
\end{eqnarray}
Since $-\Delta_\Omega$ is a positive self-adjoint operator and $1/{\ell^2}>0$, we have $-\Delta_\Omega+1/{\ell^2}$ is invertible.
It follows that
$$u_{n,1}:=\big(-\Delta_\Omega+1/{\ell^2}\big)^{-1}\cdot\big[f_{n,1}+(2 u_{n,0}-u_{n,-1})/{\ell^2}\big]$$ is the unique solution of \eqref{e:BVP_1}.
Moreover, $u_{n,1}\in C_0(\Omega)$. Successively, for $i=2,\ldots,n$, we consider
\begin{equation}\label{e:BVP_i}
(u-2u_{n,i-1}+u_{n,i-2})/{\ell^2}-\Delta_\Omega u=f_{n,i}.
\end{equation}
Similarly, the unique solution of \eqref{e:BVP_i} is
$$u_{n,i}:=\big(-\Delta_\Omega+1/{\ell^2}\big)^{-1}\cdot\big[f_{n,i}+(2 u_{n,i-1}-u_{n,i-2})/{\ell^2}\big],$$
and $u_{n,i}\in C_0(\Omega)$.

Let $u_{n,i}(x)$ be the approximation of a solution $u(t,x)$ for the problem \eqref{e:wave_eq_on_graph} at $t=t_i$.
We denote
$$\delta u_{n,i}=(u_{n,i}-u_{n,i-1})/{\ell} \qquad\mbox{ for } i=0,1,\ldots, n,$$ and
$$\delta^2 u_{n,i}=(\delta u_{n,i}-\delta u_{n,i-1})/{\ell} \qquad\mbox{ for }i=1,\ldots, n.$$
Then we can replace $\partial_t u(t_i,x)$ and $\partial_t^2 u(t_i,x)$ by $\delta u_{n,i}(x)$ and $\delta^2 u_{n,i}(x)$, respectively. It follows from \eqref{e:BVP_1} and \eqref{e:BVP_i} that
\begin{equation}\label{e:equation_on_minimum}
\delta^2 u_{n,i}-\Delta_\Omega u_{n,i}=f_{n,i}\qquad\mbox{ for }i=1,\ldots, n.
\end{equation}
Define Rothe's functions $u^{(n)}(t,x)$ from $[0,T]$ to $C_0(\Omega)$ by
\begin{equation}\label{e:defi_u(n)}
u^{(n)}(t,x)=u_{n,i-1}(x)+(t-t_{i-1})\cdot\delta u_{n,i}(x)\quad t\in[t_{i-1}, t_i],\ i=1,\ldots,n.
\end{equation}
For $i=1,\ldots, n$, let $Q_{T,i}:=[t_{i-1}, t_i]\times\Omega$ and $\widetilde{Q}_{T,i}:=(t_{i-1}, t_i]\times\Omega$.
Define
\begin{equation*}\label{e:defi_delta_u(n)}
\delta u^{(n)}(t,x)=\delta u_{n,i-1}(x)+(t-t_{i-1})\cdot\delta^2 u_{n,i}(x)
\quad\mbox{for }(t,x)\in Q_{T,i},
\end{equation*}
and some step functions
\begin{equation*}\label{e:defi_delta_overline_u(n)}
\overline{u}^{(n)}(t,x)=
\left\{
\begin{aligned}
&u_{n,i}(x),\quad&&\mbox{if }(t,x)\in\widetilde{Q}_{T,i},\\
&g(x),\quad&&\mbox{if }(t,x)\in[-\ell,0]\times\Omega^\circ,\\
&0,\quad&&\mbox{if }(t,x)\in[-\ell,0]\times\partial\Omega,
\end{aligned}
\right.
\end{equation*}
\begin{equation*}\label{e:defi_delta_overline_delta_u(n)}
\delta \overline{u}^{(n)}(t,x)=
\left\{
\begin{aligned}
&\delta u_{n,i}(x),\quad&&\mbox{if }(t,x)\in\widetilde{Q}_{T,i},\\
&h(x),\quad&&\mbox{if }(t,x)\in[-\ell,0]\times\Omega^\circ,\\
&0,\quad&&\mbox{if }(t,x)\in[-\ell,0]\times\partial\Omega,
\end{aligned}
\right.
\end{equation*}
\begin{equation*}\label{e:defi_f(n)}
f^{(n)}(t,x)=
\left\{
\begin{aligned}
&f_{n,i}(x),\quad&&\mbox{if }(t,x)\in\widetilde{Q}_{T,i},\\
&f(0,x),\quad&&\mbox{if }(t,x)\in\{0\}\times\Omega^\circ,\\
&0,\quad&&\mbox{if }(t,x)\in\{0\}\times\partial\Omega.
\end{aligned}
\right.
\end{equation*}

In order to prove that Rothe's function $u^{(n)}(t,x)$ converges to a solution $u(t,x)$ of \eqref{e:wave_eq_on_graph}, we need to give some priori estimates. In the rest of this subsection, we assume that \eqref{eq:condition_on_f} holds.
\begin{lemma}\label{L:estimates_on_delta_ui_and_so_on}
There exist an integer $N'>0$ and some positive constants $C_1, C_2, C_3$,
depending only on $\Omega$ and $T$, such that the following inequalities hold for all $n\ge N'$ and $i=1,\ldots, n$,
\begin{eqnarray*}\label{e:bounded_of_delta ui_and_ui_and_delta2_ui}
\|\delta u_{n,i}\|^2_{L^2(\Omega)}+\|\nabla u_{n,i}\|^2_{L^2(\Omega)}\le C_1,
\quad\|u_{n,i}\|^2_{L^2(\Omega)}\le C_2, \quad\|\delta^2 u_{n,i}\|^2_{L^2(\Omega)}\le C_3.
\end{eqnarray*}
\end{lemma}

\begin{proof}
Let $c':=\|f(0,\cdot)\|^2_{L^2(\Omega^\circ)}$.
Then \eqref{eq:condition_on_f} implies that
\begin{eqnarray*}
\|f(t,\cdot)\|^2_{L^2(\Omega^\circ)}\le cT^{2\alpha}+c' \qquad\mbox{ for any }t\in[0,T].
\end{eqnarray*}

By \eqref{e:equation_on_minimum}, we get that for $i=1,\ldots,n$ and any $v\in W_0^{1,2}(\Omega)$,
\begin{equation}\label{e:delta2_ui_is_a_uni_sol}
\int_{\Omega^\circ}\delta^2u_{n,i}\cdot v\,d\mu-\int_{\Omega^\circ}\Delta_\Omega u_{n,i}\cdot v\,d\mu=\int_{\Omega^\circ}f_{n,i}\cdot v\,d\mu.
\end{equation}
In \eqref{e:delta2_ui_is_a_uni_sol}, letting $v=\delta u_{n,i}$ and using Lemma \ref{L:Green_formula}, we get
\begin{eqnarray*}
\begin{aligned}
&\|\nabla u_{n,i}\|^2_{L^2(\Omega)}+(1-\ell)\|\delta u_{n,i}\|^2_{L^2(\Omega)}\\
\le& \|\nabla u_{n,i-1}\|^2_{L^2(\Omega)}+\|\delta u_{n,i-1}\|^2_{L^2(\Omega)}+\ell\|f_{n,i}\|^2_{L^2(\Omega^\circ)}.
\end{aligned}
\end{eqnarray*}
Choosing a positive integer $N'$ such that $\ell< 1$ holds for any $n\ge N'$, we obtain
\begin{eqnarray*}
\begin{aligned}
&(1-\ell)^i\big(\|\nabla u_{n,i}\|^2_{L^2(\Omega)}+\|\delta u_{n,i}\|^2_{L^2(\Omega)}\big)\\
\le&\|\nabla u_{n,0}\|^2_{L^2(\Omega)}+\|\delta u_{n,0}\|^2_{L^2(\Omega)}+\sum_{k=1}^i\ell(1-\ell)^{k-1}\|f_{n,k}\|^2_{L^2(\Omega^\circ)},
\end{aligned}
\end{eqnarray*}
and so
\begin{eqnarray*}
\begin{aligned}
&\|\nabla u_{n,i}\|^2_{L^2(\Omega)}+\|\delta u_{n,i}\|^2_{L^2(\Omega)}\\
\le&(1-\ell)^{-n}\big(\|\nabla u_{n,0}\|^2_{L^2(\Omega)}+\|\delta u_{n,0}\|^2_{L^2(\Omega)}\big)
   +\ell(1-\ell)^{-n}\sum_{k=1}^i\|f_{n,k}\|^2_{L^2(\Omega^\circ)}\\
\le&e^T \big(\|\nabla u_{n,0}\|^2_{L^2(\Omega)}+\|\delta u_{n,0}\|^2_{L^2(\Omega)}\big)+(cT^{2\alpha}+c')Te^T=:C_1,
\end{aligned}
\end{eqnarray*}
where we use the fact that
\begin{eqnarray*}
\lim_{n\to \infty}\frac{1}{(1-\ell)^{n}}=\lim_{n\to\infty}\Big(1+\frac{1}{-\frac{n}{T}}\Big)^{-\frac{n}{T}\cdot T}=e^T.
\end{eqnarray*}

It is easy to see that for any $i\in\{1,\ldots,n\}$, $u_{n,i}$ is bounded on $\Omega$. In fact,
by $\|\nabla u_{n,i}\|^2_{L^2(\Omega)}\le C_1$, we get
\begin{eqnarray*}
|u_{n,i}(y)-u_{n,i}(x)|\le \sqrt{2C_1/\omega_{min}}\qquad\mbox{ for any }x\in \Omega \mbox{ and any } y\sim x,
\end{eqnarray*}
where $\omega_{min}:=\min_{x,y\in \Omega}\omega_{xy}$. On the other hand, for any $x\in\Omega$ and some $x'\in\partial\Omega$, we can choose a shortest path on $G$ from $x$ to $x'$:
$$x=x_1\sim x_2\cdots\sim x_{k-1}\sim x_k=x'.$$
Using the triangle inequality and $u_{n,i}(x')=0$, we get
\begin{eqnarray*}
|u_{n,i}(x)|
\le |u_{n,i}(x_1)-u_{n,i}(x_2)|+\cdots+|u_{n,i}(x_{k-1})-u_{n,i}(x_k)|
\le k\sqrt{2C_1/\omega_{min}}.
\end{eqnarray*}
This proves that $u_{n,i}$ is bounded on $\Omega$,
and hence
\begin{equation*}
\|u_{n,i}\|^2_{L^2(\Omega)}\le C_2.
\end{equation*}

Since $|\Delta_\Omega u_{n,i}|^2\le D_\mu|\nabla u_{n,i}|^2$, we get $\|\Delta_\Omega u_{n,i}\|^2_{L^2(\Omega)}\le C_1 D_\mu$,
which, together with \eqref{e:equation_on_minimum}, yields
\begin{eqnarray*}\label{e:bounded_of_delta2 ui}
\|\delta^2 u_{n,i}\|^2_{L^2(\Omega)}\le C_3.
\end{eqnarray*}
This completes the proof.
\end{proof}

By Lemma \ref{L:estimates_on_delta_ui_and_so_on}, we can get the following estimates.
\begin{lemma}
For any $t\in[0,T]$ and any $n\ge N'$, there exist some positive constants $C_4, C_5$ depending only on $\Omega$ and $T$, such that
\begin{eqnarray}\label{e:bounded_1}
\begin{aligned}
&\|u^{(n)}(t,\cdot)\|_{L^{2}(\Omega)}+\|\overline{u}^{(n)}(t,\cdot)\|_{L^{2}(\Omega)}+\|\delta u^{(n)}(t,\cdot)\|_{L^2(\Omega)}\\
+&\|\delta \overline{u}^{(n)}(t,\cdot)\|_{L^2(\Omega)}+\|\partial_t\big(\delta u^{(n)}\big)(t,\cdot)\|_{L^2(\Omega)}\le C_4,
\end{aligned}
\end{eqnarray}
and
\begin{eqnarray}\label{e:bounded_2}
\|u^{(n)}(t,\cdot)-\overline{u}^{(n)}(t,\cdot)\|_{L^2(\Omega)}
+\|\delta u^{(n)}(t,\cdot)-\delta \overline{u}^{(n)}(t,\cdot)\|_{L^2(\Omega)}\le C_5/n.
\end{eqnarray}
\end{lemma}

\begin{proof}
Using Lemma \ref{L:estimates_on_delta_ui_and_so_on}, it is easy to see that \eqref{e:bounded_1} holds.
Now, we prove \eqref{e:bounded_2}. If $t=t_{i-1}$, then
$\|u^{(n)}(t,\cdot)-\overline{u}^{(n)}(t,\cdot)\|_{L^2(\Omega)}
+\|\delta u^{(n)}(t,\cdot)-\delta\overline{u}^{(n)}(t,\cdot)\|_{L^2(\Omega)}=0$.
If $t\in(t_{i-1},t_i]$, then
\begin{eqnarray*}
\begin{aligned}
&\|u^{(n)}(t,\cdot)-\overline{u}^{(n)}(t,\cdot)\|_{L^2(\Omega)}
+\|\delta u^{(n)}(t,\cdot)-\delta \overline{u}^{(n)}(t,\cdot)\|_{L^2(\Omega)}\\
\le&2\ell\cdot\big(\|\delta u_{n,i}\|_{L^2(\Omega)}+\|\delta^2 u_{n,i}\|_{L^2(\Omega)}\big)\\
\le&\frac{2T(\sqrt{C_1}+\sqrt{C_3})}{n}.
\end{aligned}
\end{eqnarray*}
This completes the proof of the results.
\end{proof}

\begin{lemma}\label{L:weak_convergence}
There exist a function $u\in L^2(\Omega)$ satisfying $\partial_tu, \partial_t^2u\in L^2(\Omega)$, and two subsequences $\{u^{(n)}\}$, $\{\overline{u}^{(n)}\}$, denoted by themselves, such that
\begin{enumerate}
\item[(a)] $u^{(n)}\to u$ and $\overline{u}^{(n)}\to u$ on $[0,T]\times\Omega$;
\item[(b)] $\delta u^{(n)}\to \partial_tu$ and $\delta \overline{u}^{(n)}\to \partial_tu$ on $[0,T]\times\Omega$;
\item[(c)] $\partial_t(\delta u^{(n)})\to \partial_t^2u$ on $[0,T]\times\Omega$.
\end{enumerate}
\end{lemma}

\begin{proof}
(a) According to \eqref{e:bounded_1}, there exist two subsequences $\{u^{(n)}\}$, $\{\overline{u}^{(n)}\}$, and two functions $u,$ $\overline{u}$ satisfying for any $t\in[0,T]$,
\begin{eqnarray*}
u^{(n)}(t,\cdot)\to u(t,\cdot)\quad\mbox{and}\quad
\overline{u}^{(n)}(t,\cdot)\to \overline{u}(t,\cdot)\quad\mbox{ in }L^{2}(\Omega).
\end{eqnarray*}
This leads to
\begin{eqnarray*}
u^{(n)}(t,x)\to u(t,x)\quad\mbox{and}\quad
\overline{u}^{(n)}(t,x)\to \overline{u}(t,x)\quad\mbox{ for }(t,x)\in[0,T]\times\Omega.
\end{eqnarray*}
Combining this with \eqref{e:bounded_2}, we get for any $t\in[0,T]$,
\begin{eqnarray*}
 \|u(t,\cdot)-\overline{u}(t,\cdot)\|^2_{L^2(\Omega)}
=\lim_{k\to\infty}\|u^{(n)}(t,\cdot)-\overline{u}^{(n)}(t,\cdot)\|^2_{L^2(\Omega)}
=0,
\end{eqnarray*}
and so $u=\overline{u}$. Moreover, for any $t\in [0,T]$, $u(t,x)=0$ on $\partial\Omega$.

(b) Similar to (a), we can extract two subsequences $\{\delta u^{(n)}\}$ and $\{\delta \overline{u}^{(n)}\}$ such that
for some function $w\in L^2(\Omega)$,
\begin{eqnarray*}
\delta u^{(n)}(t,x)\to w(t,x)\quad\mbox{and}\quad\delta \overline{u}^{(n)}(t,x)\to w(t,x)
\qquad\mbox{ on }[0,T]\times\Omega.
\end{eqnarray*}

Next, we prove that $w=\partial_tu$. For $i=1,\ldots, n$ and any $t\in[t_{i-1},t_i]\subseteq[0,T]$,
\begin{eqnarray}\label{e:deri_ut}
u^{(n)}(t,x)-g(x)
&=&\int_0^{t_1}\partial_su^{(n)}(s,\cdot)\,ds+\cdots+\int_{t_{i-2}}^{t_{i-1}}\partial_su^{(n)}(s,\cdot)\,ds\notag\\
&~&+\int_{t_{i-1}}^t\partial_su^{(n)}(s,\cdot)\,ds\notag\\
&=&\int_0^{t_1}\delta u^1(\cdot)\,ds+\cdots+\int_{t_{i-2}}^{t_{i-1}}\delta u_{n,i-1}(\cdot)\,ds+\int_{t_{i-1}}^t\delta u_{n,i}(\cdot)\,ds\notag\\
&=&\int_0^t\delta \overline{u}^{(n)}(s,\cdot)\,ds.
\end{eqnarray}
Since
$$|\delta \overline{u}^{(n)}(t,x)|^2\mu_0\le\|\delta \overline{u}^{(n)}(t,\cdot)\|_{L^2(\Omega)}\le C_4\quad\mbox{ for }(t,x)\in[0,T]\times \Omega,$$
we get
\begin{eqnarray*}
|\delta \overline{u}^{(n)}(t,x)|\le \sqrt{C_4 /{\mu_0}}\qquad\mbox{ on }[0,T]\times \Omega,
\end{eqnarray*}
where $\mu_0=\min_{x\in \Omega}\mu(x)$. Using Dominated Convergence Theorem, we get
\begin{eqnarray}\label{e:fact_delta_over_u_conver_to_w_0t}
\int_0^t\delta \overline{u}^{(n)}(s,\cdot)\,ds\to \int_0^t w(s,\cdot)\,ds.
\end{eqnarray}
Taking the limit as $n\to\infty$ in \eqref{e:deri_ut}, and using (a), we get
\begin{equation*}
u(t,x)-g(x)=\int_0^t w(s,\cdot)\,ds.
\end{equation*}
Hence $w=\partial_t u$ and $u(0,x)=g(x)$.

(c) Similar to (b), there exits a subsequence $\{\partial_t(\delta u^{(n)})\}$ satisfying
\begin{equation*}
\partial_t (\delta u^{(n)}) \to \partial_t^2u\qquad\mbox{on }[0,T]\times\Omega.
\end{equation*}
Moreover, $\partial_t u(0,x)=h(x)$. In the proof, we use the fact that for $t\in[0,T]$,
\begin{eqnarray}\label{e:fact_partial_delta_u_conver_to_partial_tt_u_0t}
\int_0^t \partial_s(\delta u^{(n)})(s,x)\,ds \to \int_0^t \partial^2_s u(s,x)\,ds\quad\mbox{on }\Omega,
\end{eqnarray}
whose proof is the same as that of \eqref{e:fact_delta_over_u_conver_to_w_0t}. This completes the proof.
\end{proof}

\begin{lemma}\label{L:weak_convergence_on_0_T}
For any $n\ge N'$, the following results hold:
\begin{enumerate}
\item[(a)] $\int_0^T \Delta_\Omega \overline{u}^{(n)}(t,x)\,dt\to \int_0^T \Delta_\Omega u(t,x)\,dt$ on $\Omega^\circ$;
\item[(b)] $\int_0^T f^{(n)}(t,x)\,dt\to \int_0^T f(t,x)\,dt$ on $\Omega^\circ$.
\end{enumerate}
\end{lemma}
\begin{proof}
(a) From Lemma \ref{L:weak_convergence}(a), we get $\Delta_\Omega \overline{u}^{(n)}\to \Delta_\Omega u$ on $[0,T]\times\Omega^\circ$. Combining this with \eqref{e:bounded_1} and Dominated Convergence Theorem, the result holds.

(b) The proof is the same as that of (a).
\end{proof}

\subsection{Proof of Theorem \ref{T:sol_of_wave_equ}}

In this subsection, we use the notation defined in Subsection \ref{SS:estimates}.
Under the derivation of Subsection \ref{SS:estimates}, we prove that \eqref{e:wave_eq_on_graph} has a unique solution.
\begin{proof}[Proof of Theorem~\ref{T:sol_of_wave_equ}]
\noindent\textit{Existence.} 
By \eqref{e:equation_on_minimum}, we get
\begin{eqnarray}\label{e:int_of_all_functions}
\int_0^T\partial_t(\delta u^{(n)})(t,x)\,dt-\int_0^T\Delta_\Omega \overline{u}^{(n)}(t,x)\,dt=\int_0^T f^{(n)}(t,x)\,dt
\quad\mbox{on }\Omega^\circ.
\end{eqnarray}
Let $u$ be the limit function in Lemma \ref{L:weak_convergence}.
Taking the limit as $n\to\infty$ in \eqref{e:int_of_all_functions}, using \eqref{e:fact_partial_delta_u_conver_to_partial_tt_u_0t} and Lemma \ref{L:weak_convergence_on_0_T}, we get
\begin{eqnarray*}
\int_0^T\big(\partial^2_t u(t,x)-\Delta_\Omega u(t,x)-f(t,x)\big)\,dt=0\quad\mbox{on }\Omega^\circ.
\end{eqnarray*}
From Lemma \ref{L:weak_convergence}, we get $u(0,x)=g(x)$, $\partial_tu(0,x)=h(x)$, and
$u(t,x)=0$ on $[0,T]\times\partial\Omega$.
Since $T>0$ is arbitrary, $u$ is a solution of \eqref{e:wave_eq_on_graph}.

\noindent\textit{Uniqueness.}
If $u_1$ and $u_2$ both satisfy \eqref{e:wave_eq_on_graph}, then $w:=u_1-u_2$ satisfies
\begin{equation}\label{e:uniqueness_of_wave_eq}
\left\{
\begin{aligned}
&\partial_t^2w-\Delta_\Omega w=0, \quad&(t,x)\in[0,\infty)\times \Omega^\circ,\\
&w(0,x)=0,\quad&x\in\Omega^\circ,\\
&\partial_tw(0,x)=0,\quad&x\in\Omega^\circ,\\
&w(t,x)=0,\quad&(t,x)\in[0,\infty)\times\partial\Omega.
\end{aligned}
\right.
\end{equation}
Let
\begin{eqnarray*}\label{e:e(t)}
e(t)=\int_\Omega|\nabla w(t,x)|^2\,d\mu+\int_{\Omega^\circ}|\partial_tw(t,x)|^2\,d\mu\qquad\mbox{ for }t\in[0,\infty).
\end{eqnarray*}
Then $e(0)=0$ and
\begin{eqnarray*}
\begin{aligned}
e'_+(0)
&=\lim_{t\ra 0^+}\frac{\int_\Omega|\nabla w(t,x)|^2\,d\mu+\int_{\Omega^\circ}|\partial_tw(t,x)|^2\,d\mu}{t}\\
&=2\int_{\Omega}\Big(\lim_{t\to 0^+}\Gamma(w(t,x), \partial_tw(t,x))\Big)\,d\mu
+2\int_{\Omega^\circ}\Big(\lim_{t\ra 0^+} \partial_tw(t,x)\cdot \partial_t^2w(t,x)\Big)\,d\mu\\
&=0.
\end{aligned}
\end{eqnarray*}
Moreover, for $t\in(0,\infty)$,
\begin{eqnarray*}
\begin{aligned}
e'(t)
&=2\int_\Omega \Gamma(w,\partial_tw)\,d\mu +2\int_{\Omega^\circ}\partial_tw\cdot \partial_t^2w\,d\mu\\
&=-2\int_{\Omega^\circ} \Delta_\Omega w\cdot \partial_tw\,d\mu+2\int_{\Omega^\circ}\partial_tw\cdot \partial_t^2w\,d\mu\\
&=0.
\end{aligned}
\end{eqnarray*}
Thus $e(t)\equiv 0$ for $t\in[0,\infty)$, and so
$$|\nabla w|(t,x)\equiv0\mbox{ for }(t,x)\in[0,\infty)\times\Omega\quad\mbox{and}\quad
\partial_tw(t,x)\equiv0\mbox{ for }(t,x)\in[0,\infty)\times\Omega^\circ.$$
For a fixed $t\in[0,\infty)$, the facts $|\nabla w|(t,x)\equiv 0$ and $\Omega$ is connected imply that $w(t,x)\equiv \mbox{constant}$ for any $x\in \Omega$. Fixed $x\in \Omega^\circ$, it follows from $\partial_tw(t,x)\equiv 0$ that $w(t,x)\equiv \mbox{constant}$ for any $t\in[0,\infty)$. Combining these with \eqref{e:uniqueness_of_wave_eq}, we get $w(t,x)\equiv0$ for $(t,x)\in[0,\infty)\times\Omega$.
Applying this argument to $w=u_1-u_2$ gives that $u_1=u_2$.
\end{proof}

\section{Proof of Theorem~\ref{T:sol_of_nonlinear_wave_equ}}\label{S:proof_of_main_results_nonlinear}

\begin{proof}[Proof of Theorem \ref{T:sol_of_nonlinear_wave_equ}]
By the method of variation of constant, we get
\begin{eqnarray*}
\begin{aligned}
u(t,x)
=&\sum_{k=1}^N\frac{1}{\sqrt{\lambda_k}}\int_0^t\sin(\sqrt{\lambda_k}(t-s))b_k(s)\,ds\cdot \varphi_k(x)\\
&+\sum_{k=1}^Ng_k\cos(\sqrt{\lambda_k}t)\varphi_k(x)+\sum_{k=1}^N\frac{1}{\sqrt{\lambda_k}}h_k\sin(\sqrt{\lambda_k}t)\varphi_k(x).
\end{aligned}
\end{eqnarray*}
It is obvious that $u(t,x)$ satisfies \eqref{e:wave_eq_on_graph}. It follows from Theorem \ref{T:sol_of_wave_equ} that \eqref{e:wave_eq_on_graph} has a unique solution. This completes the proof of Theorem \ref{T:sol_of_nonlinear_wave_equ}.
\end{proof}


\textbf{Acknowledgement}
This research was supported by the National Science Foundation of China [grant 12071245].

\section*{References}
\bibliographystyle{amsalpha}

\end{document}